\pdfoutput=1
\documentclass{article}
\usepackage{amsmath,amsfonts,amsthm}
\usepackage{mathrsfs}
\usepackage{graphicx}
\usepackage{hyperref}
\usepackage{tabu}
\usepackage{booktabs}
\usepackage{makecell}

\newtheorem{theorem}{Theorem}[section]

\newtheorem{definition}[theorem]{Definition}
\theoremstyle{remark}

\newcommand{\R}{\mathbb{R}}
\newcommand{\Fil}{\mathcal{F}}
\newcommand{\Prob}{\mathbb{P}}

\newcommand{\E}{\mathbb{E}}
\newcommand{\F}{\mathscr{F}}
\newcommand{\G}{\mathscr{G}}

\newcommand\mgape[1]{\gape{$\vcenter{\hbox{#1}}$}}

\numberwithin{equation}{section}

\title{Parareal algorithms applied to stochastic differential equations with conserved quantities}
\author{Liying Zhang\footnote{School of Mathematical Science, China University of Mining and Technology, Beijing 100083, China. (lyzhang@lsec.cc.ac.cn)}, Weien Zhou\footnote{College of Science, National University of Defense Technology, Changsha 410073, China. (weienzhou@nudt.edu.cn)} and Lihai Ji\footnote{Institute of Applied Physics and Computational Mathematics, Beijing, China. (jilihai@lsec.cc.ac.cn)}}
\date{}

\begin{document}
\maketitle
\begin{abstract}
	In this papers, we couple the parareal algorithm  with projection methods of the trajectory on a specific manifold, defined by the preservation of some conserved quantities of the  differential equations. First, projection methods are introduced as the coarse and fine propagators. Second, we also apply the projection methods for systems with conserved quantities in the correction step of original parareal algorithm. Finally, three numerical experiments are performed by different kinds of algorithms to show the property of convergence in iteration, and preservation in conserved quantities of model systems.

\end{abstract}

\section{Introduction}\label{s;intro}
Designing high efficient algorithms is an important subject of numerical computation due to the computational time
and memory issues in the solution
of large scale problems. 
The technique of parallel algorithms was paid more and more attention in past few years, containing domain decomposition method in spatial direction and the parallel in time direction generally.
The parareal algorithm, our focus in the sequel, was first introduced by Lions, Maday, and Turinici in 2001 \cite{Lions2001}, further work modified by Bal and Maday in \cite{Bal2002}, and has attracted vast attention in the last decade.
Compared with other parallel approaches, this algorithm belongs to time-parallel category. 
The general idea of parareal algorithm contains roughly  three steps  as follows.
First, we obtain an approximate solution on a coarse time-step by a rough solver. Second,  we use another more accurate solver to get the approximation on each coarse time interval (splitting the coarse time interval into more fine time domain) performed by parallel with initial values computed in the first step. Finally, combining the values of the above two steps in the coarse time grids, we obtain a new approximation value by a prediction and correction iteration.
In general, this algorithm has higher parallel performance and is easy to perform, which motivates the development of efficient parallel methods for time  dependent problems. Since the parareal algorithm was proposed, many efforts have been made to analyze it theoretically \cite{Wu2015} and numerically, which verify the effectiveness of the parareal algorithm for a large various of problems, including control theory\cite{Maday2002}, Navier-Stokes problem\cite{2005Barth} and Hamiltonian differential equations\cite{Dai2013,Gander2014} for instance.

Stochastic differential equations (SDEs) have attracted considerable attention in order to obtain much more realistic mathematical models in many scientific disciplines, such as physics, molecular biology, population dynamics and finance \cite{Gardiner2009,Oksendal2003}.
However,
it is difficult to find explicit solutions of SDEs analytically; therefore, there has been tremendous interest in developing effective and reliable numerical methods for SDEs  (e.g. \cite{Burrage2004,Higham2001,Kloeden1992} and references therein).
It is also a significant issue whether some geometric features of SDEs are preserved in performing reliable numerical methods, especially for long-time simulation, which is as important as the deterministic case \cite{Hairer2006,Milstein2004}. 
In practice, they are time consuming, so the parallel techniques can be considered to speed up the original integrator.
For stochastic problem, the application of parallel algorithm  are relatively few. For example, the parareal algorithm has been applied to stochastic ordinary differential equations with filter problems \cite{Bal2003} and stochastic models in chemical kinetics \cite{Engblom2009}.
However, to the best of our knowledge, no results on parareal algorithm focusing on stochastic differential equations with conserved quantities.
In order to apply the parareal algorithm to SDEs with conserved quantities, as mentioned in \cite{Dai2013,Dai2013a,Gander2014},
the original algorithm are not able to share this kind of conservative property, namely, the preservation of conserved quantities along the sample path of the exact solution, even though when the coarse and fine integrators all have adequate conservative properties. Therefore, the behavior of long time numerical simulation is not enjoyed as the original system itself has. In this paper, we mainly utilize the projection methods for SDEs with conserved quantities as
the basic propagators and the parareal algorithm with a projection corrector, which preserve some conserved quantities of the exact flow as proposed in \cite{Dai2013}.

The rest of the paper is organized as follows. Section \ref{s;pr} briefly recalls the parareal algorithm for general time-dependent problem. Section \ref{s;pro} discusses the procedure projection methods for SDEs with conserved quantities, and gives the corresponding mean-square convergence. Next in Section \ref{s;prq}, we consider the parareal algorithm focusing on the SDEs with certain conserved quantities, which combines the ideas of the previous two sections. Finally, three typical SDE examples are chosen to perform numerical tests in Section \ref{s;numer}.

\section{The original parareal algorithm}\label{s;pr}
In this section, we first review the original parareal algorithm for a general initial-value problem:
\begin{equation}\label{2.2.1}
\left\{
\begin{aligned}
u'(t) &= f(t,u(t)), \quad t\in [0,T], \\
u(0) &= u_0,
\end{aligned}
\right.
\end{equation}
where $f: \R \times \R^d \rightarrow \R^d $ is a suitable function to ensure the well-posedness of (\ref{2.2.1}). To perform the parareal algorithm, we first divide time interval $[0,T]$ into $N$ uniform large time intervals $[T_n, T_{n+1}]$, with step-size $\Delta T = T_{n+1} - T_n$ $n=0,1,\dots,N-1$, and $N=\frac{T}{\Delta T}$. Then, we further divide every large interval $[T_n, T_{n+1}]$ into $J$ small time intervals $[t_{n+\frac{j}{J}}, t_{n+\frac{j+1}{J}}]$, $j = 0,1,\dots, J-1$. With that, two numerical propagators, the coarse propagator $\G$ and the fine propagator $\F$, are needed here. In fact, $\G$ is usually easy to solve with low convergence order and $\F$ is of high order but more expensive to compute. The parareal algorithm can be described as following.
\begin{itemize}
	\item Initialization: use the coarse propagator $\mathscr{G}$ and time-step $\Delta T$ to compute initial value $\{u_{n}^{0}\}_{n=0}^N$ sequentially
	\begin{equation*}
	\left\{
	\begin{aligned}
	u_{n+1}^{0} &= \G(T_n, u_n^{(0)},\Delta T), \quad n=0,1,\dots,N-1,\\
	u_0^{(0)} &= u_0.
	\end{aligned}\right.
	\end{equation*}
	\item For $k=0,1,\dots$
	\begin{enumerate}
		\item use the fine propagator $\F$ and small time-step $\frac{\Delta T}{J}$  to compute $\hat{u}_n$ on each sub-interval $[T_n, T_{n+1}]$ independently,  thus possibly in parallel
		\begin{equation*}
		\left\{
		\begin{aligned}
		\hat{u}_{n+\frac{j+1}{J}} &= \F(t_{n+\frac{j}{J}}, \hat{u}_{n+\frac{j}{J}},\frac{\Delta  T}{J}), \quad j=0,1,\dots,J-1,\\
		\hat{u}_n &= u_0^{(k)}.
		\end{aligned}\right.
		\end{equation*}
		\item perform sequential corrections
		\begin{equation*}
		\left\{
		\begin{aligned}
		u_{n+1}^{(k+1)} &= \G(T_n, u_n^{(k+1)},\Delta T) + \hat{u}_{n+1} - \G(T_n, u_n^{(k)},\Delta T), \quad n=0,1,\dots,N-1,\\
		u_0^{(0)} &= u_0.
		\end{aligned}\right.
		\end{equation*}
		\item If $\{u_n^{k+1}\}_{n=1}^N$ satisfies the stopping criterion, break the iteration; otherwise continue the iteration again.
	\end{enumerate}
\end{itemize}
Note that the parareal algorithm can be expressed compactly as follows:
\begin{equation}\label{e;pc}
u_{n+1}^{(k+1)} = \G(T_n, u_n^{(k+1)},\Delta T) + \F^J(T_n, u_n^{(k)},\frac{\Delta T}{J}) - \G(T_n, u_n^{(k)},\Delta T)
\end{equation}
where $\F^J$ means computing the value of $\F$ for $J$ times sequentially. It is known that $u_n^{(k)} \rightarrow u_n^*, n=0,1,\dots,N$, as $k \rightarrow + \infty$ when iteration \eqref{e;pc} converges, where $u_n^*$ is actually the result computed by the fine propagator $\F$ with small step-size $\Delta T/J$ \cite{Wu2011}. Thus, the convergence accuracy of this iterative algorithm after certain iterations is comparable to that of the fine propagator $\F$ with the small step-size $\Delta T/J$ \cite{Bal2002}. In other words, the parareal algorithm can approach to the accuracy of the fine propagator, and the computational cost only is same as the coarse propagator.

\section{Projection methods for SDEs with conserved quantities}\label{s;pro}

Consider the initial value problem for the general $d$-dimensional autonomous stochastic differential equation (SDE) in the sense of Stratonovich:
\begin{equation}\label{e;sde}
\left\{
\begin{aligned}
dX(t) &= f\big(X(t)\big)dt + \sum_{r=1}^m g_r\big(X(t)\big) \circ dW_r(t), \quad t\in [0, T],\\
X(0) &= X_0,
\end{aligned}\right.
\end{equation}
where $X(t)$ is $d$-dimensional column-vector, $W_r(t), r = 1,\dots,m $, are $m$ independent one-dimensional standard Wiener processes defined on a complete filtered probability space $(\Omega , \Fil, \Prob, \{\Fil_t\}_{t \geq 0})$ fulfilling the usual conditions, $f$ and $g_r$ are $\R^d$-valued functions satisfying the conditions under which \eqref{e;sde} has a unique solution. $X_0$ is $\Fil_{0}$-measurable random variable with $\E|X_0|^2<\infty$.

\begin{definition}\label{d:multiple}
	System \eqref{e;sde} possesses $l$ ($l\ge 1$) independent conserved quantities $I^i(x)$, $i=1,\dots ,l$, if
	\begin{equation*}
	\big(\nabla I^i(x)\big)^T f(x) = 0 \quad\text{and}\quad\big(\nabla I^i(x)\big)^T g_r(x) = 0,
	\quad r=1,\dots,m; \quad i=1,\dots,l.
	\end{equation*}
\end{definition}
If we define vector $\mathbf{I}(x) := \big(I^1(x),\dots,I^l(x)\big)^T$, then
\[\mathbf{I}'(x)f(x)=\mathbf{I}'(x)g_r(x)=\mathbf{0}, \quad r=1,\dots,m,\]
where $\mathbf{I}'(x)$ is the Jacobian matrix of $I(x)$.
If system \eqref{e;sde} possesses $l$ conserved quantities $I^i(x)$, $i=1,\dots ,l$, then by It\^{o}'s formula we have
\begin{equation*}\label{e;dI}
\begin{aligned}
dI^i\big(X(t)\big) = \nabla I^i\big(X(t)\big)^T f\big(X(t)\big)dt + \sum_{r=1}^{m} \nabla I^i\big(X(t)\big)^T g_r\big(X(t)\big)\circ dW_r(t)=0.
\end{aligned}
\end{equation*}
Then
\begin{equation*}
X(t)\in \mathcal{M}_{X_0}:=\Big\{x\in\R^d \mid I^i(x)=I^i(X_0),
i=1,\dots,l\Big\} \quad t\in[0,T], \quad \text{a.s.},
\end{equation*}
which implies that the solution $X(t)$ of this system will be confined to the invariant submanifold $\mathcal{M}_{X_0}$ generated by $I^i(x)$, $i=1,\dots,l$.

\begin{figure}[ht]
	\centering
	\includegraphics[width=1\linewidth]{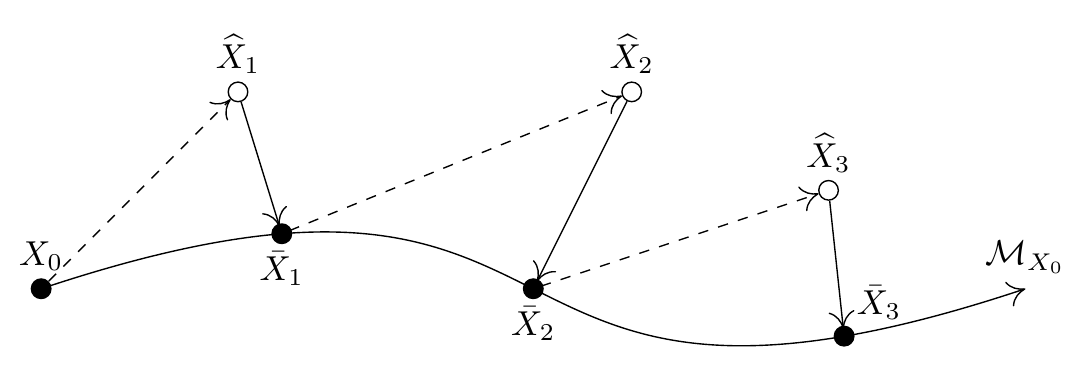}
	\caption{Basic idea of the projection methods.}
	\label{f:proj}
\end{figure}

Suppose that we have a supporting one-step method $\widehat{X}_{t,x}$, the projection method, then the process is
\begin{enumerate}
	\item Compute the one-step approximation $\widehat{X}_{t,x}$.
	\item Compute $\mathbf{\lambda}\in \R^l$ for $\bar{X}_{x,t} = \widehat{X}_{t,x} + \Phi\lambda$,\, s.t.\, $\mathbf{I}(\bar{X}_{x,t}) = \mathbf{I}(x)$.
\end{enumerate}
Here the matrix $\Phi \in \R^{d\times l}$ defines the direction of the projection, and $\mathbf{\lambda}$ is a $l$-dimensional vector chosen such that $\bar{X}_{t,x}$ belongs to the invariant manifold $\mathcal{M}_{X_0}$.  In fact $\Phi$ is not unique, and here we choose  $\Phi=\bigl(\mathbf{I}'(\widehat{X}_{t,x})\bigr)^T$, which is transpose of the Jacobian matrix of $\mathbf{I}(\cdot)$ at $\widehat{X}_{x,t}$.
The general idea of the projection methods is shown in Fig. \ref{f:proj}.

The convergence in the mean-square of this kind of projection methods is listed below.
\begin{theorem}	\cite{Zhou2016}
	Suppose that system \eqref{e;sde} possesses $l$ independent conserved quantities $I^i(x),i=1,\dots,l$.
	Also assume that a supporting method $\widehat{X}$ applying to \eqref{e;sde} satisfies
	\begin{equation}\label{c:p1}
	|\E(X_{t,x}(t+h) - \widehat{X}_{t,x}(t+h))| \leq K(1+|x|^2)^{1/2} h^{p+1},
	\end{equation}
	\begin{equation}\label{c:p2}
	\big(\E|X_{t,x}(t+h) - \widehat{X}_{t,x}(t+h)|^2\big)^{1/2} \leq K(1+|x|^2)^{1/2} h^{p+\frac{1}{2}},
	\end{equation}
	with mean-square order $p$. Assume that $\nabla I^i$ satisfies global Lipschitz condition and has uniformly bounded derivatives up to order 2, $|\nabla I^i|$ has a positive lower bound and $\big(|\nabla I^i|^2\big)^{-1}$ has bounded derivative near the invariant manifold. Then the projection method $\bar{X}$ using the supporting method $\widehat{X}$ will also have mean-square order p as well.
\end{theorem}

\section{Parareal algorithm for SDEs with conserved quantities}\label{s;prq}
For SDEs with conserved quantities, both theoretical and numerical results show that the original parareal algorithm in Section \ref{s;pr} is unable to deal with this kind of problem in long time simulation \cite{Bal2003,Dai2013a,Gander2014}, so we need other technique to deal with it. Of course, the projection method is a natural choice in order to preserve the conserved quantities of system. Even though we can choose the projection methods described in Section \ref{s;pro} as propagators $\G$ and $\F$ in the original parareal algorithm, after sequential corrections, the new iterations cannot preserve the conversed quantities any longer. Thus, what we need is another projection step to ensure that the approximations in every iteration preserve the conserved quantities as well.

To be precise, we list the corresponding parareal algorithm with projection for SDEs. As in Section \ref{s;pr}, we have the coarse propagator $\G$ and the fine propagator $\F$ for SDE \eqref{e;sde}, but here they converge in the mean-square sense.
\begin{itemize}
	\item Initialization: use the coarse propagator $\G$ and time-step $\Delta T$ to compute initial value $\{X_{n}^{0}\}_{n=0}^N$ sequentially
	\begin{equation*}
	\left\{
	\begin{aligned}
	X_{n+1}^{0} &= \G(T_n, X_n^{(0)},\Delta T), \quad n=0,1,\dots,N-1,\\
	X_0^{(0)} &= X_0.
	\end{aligned}\right.
	\end{equation*}
	\item For $k=0,1,...$
	\begin{enumerate}
		\item use the fine propagator $\F$ and small time-step $\frac{\Delta T}{J}$ to compute $\hat{X}_n$ on each sub-interval $[T_n, T_{n+1}]$ independently
		\begin{equation*}
		\left\{
		\begin{aligned}
		\hat{X}_{n+\frac{j+1}{J}} &= \F(t_{n+\frac{j}{J}}, \hat{X}_{n+\frac{j}{J}},\Delta  T/J), \quad j=0,1,\dots,J-1,\\
		\hat{X}_n &= X_0^{(k)}
		\end{aligned}\right.
		\end{equation*}
		\item perform sequential corrections
		\begin{equation*}
		\left\{
		\begin{aligned}
		X_{n+1}^{(k+1)} &= \pi_{\mathcal{M}_{X_0}} \Big(\G(T_n, X_n^{(k+1)},\Delta T) + \hat{X}_{n+1} - \G(T_n, X_n^{(k)},\Delta T)\Big), \quad n=0,1,\dots,N-1,\\
		X_0^{(0)} &= X_0,
		\end{aligned}\right.
		\end{equation*}
		where $\pi_{\mathcal{M}_{X_0}}(\cdot)$ denotes the projection operator.
		\item If $\{X_n^{k+1}\}_{n=1}^N$ satisfy the stopping criterion, break the loop; otherwise continue the iteration again.
	\end{enumerate}
\end{itemize}

Note that, in the sequential correction step, we couple an additional projection operator applied to the original parareal algorithm so that the new iteration confined on the same invariant manifold, which implies it can preserves the conserved quantities of the system. Furthermore, $X_{n}^{(k)}$ converge to $\F$ with projection $\pi_{\mathcal{M}_{X_0}}$, denoted by $\F{\pi_{\mathcal{M}_{X_0}}}$, instead of the fine propagator $\F$.

\section{Numerical experiments}\label{s;numer}
In this section, we perform several typical numerical examples by utilizing different parareal algorithms, with or without projection procedure. In order to  investigate the convergence property of these algorithms for SDEs with conserved quantities through numerical tests, we consider the following schemes:
\begin{itemize}
	\item Euler-Maruyama scheme (Euler, EulerP)
	\item Milstein scheme (Mil, MilP)
	\item Mid-point scheme (Mid, MidP)
	\item It\^{o}-Taylor order 1.5 scheme (T32, T32P)
	\item It\^{o}-Taylor order 2 scheme (T2, T2P)
\end{itemize}
where the suffix P means the projection method introduced in Section \ref{s;pro}.
That is to say, we use these schemes both for the coarse propagator $\G$ and the fine propagator $\F$, respectively.

The mean-square error is applied as  the stopping criterion of these parareal algorithms:
\begin{equation}\label{e;stopping}
(\E |X_N^{(k)} - X^*_N|^2)^{1/2} \le 10^{-12},
\end{equation}
where $X^*_N$ denotes the last step approximation computed with small step-size $\frac{\Delta T}{J}$ by the fine propagator $\F$ for original parareal, or by $\F$ with projection for parareal with projection in Section \ref{s;prq}. The expectation here is simulated by computing the average of 1000 sample paths.

\begin{table}[htbp]
	\centering
	\caption{Line styles.}
	\label{tab;linestyle}
	\begin{tabular}{@{}lcc@{}}
		\toprule
		& \multicolumn{2}{c}{Projection} \\
		\cmidrule(l){2-3}
		Style & Propagators &    Correction     \\ \midrule
		\mgape{\includegraphics{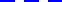}}  &  $\times$  &     $\times$      \\
		\mgape{\includegraphics{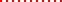}}  &  $\times$  &    \checkmark     \\
		\mgape{\includegraphics{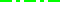}}  & \checkmark &     $\times$      \\
		\mgape{\includegraphics{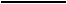}}  & \checkmark &    \checkmark     \\ \bottomrule
	\end{tabular}
\end{table}

In Table \ref{tab;linestyle}, we list four line styles to make a distinction among different algorithms in figures later. $\times$ and \checkmark denote whether the projection technique is used in propagators (both coarse and fine) and the sequential correction. For instance, if we check the type of  Euler scheme, then the solid line (the fourth style in Table \ref{tab;linestyle}) signifies that we apply the EulerP scheme in $\G$ and $\F$, and make use of the parareal algorithm with projection.

\subsection{Kubo oscillator}
Our first example is a two-dimensional linear SDEs of this form
\begin{equation}\label{kubo}
\left\{
\begin{aligned}
dX_1(t) &= -X_2(t) dt - c X_2(t) \circ dW(t),\\
dX_2(t) &= \phantom{-} X_1(t) dt + c X_1(t) \circ dW(t),
\end{aligned}\right.
\end{equation}
where $c$ is a real-valued parameter. Note that it is also called Kubo oscillator and is a typical stochastic Hamiltonian system with multiplicative noise \cite{Cohen2014,Milstein2002mul}. It is easy to check that \eqref{kubo} has a quadratic conserved quantity
\begin{equation}\label{e;I_kubo}
I(x,y) = \frac{1}{2}(x^2 + y^2),
\end{equation}
which is also its Hamiltonian function and forms a circle as the invariant submanifold in its phase space. In this example, we choose $c=0.5$, and the initial value $X(0)=(1,0)$.

\begin{figure}[htbp]
	\centering
	\includegraphics[width=0.9\textwidth]{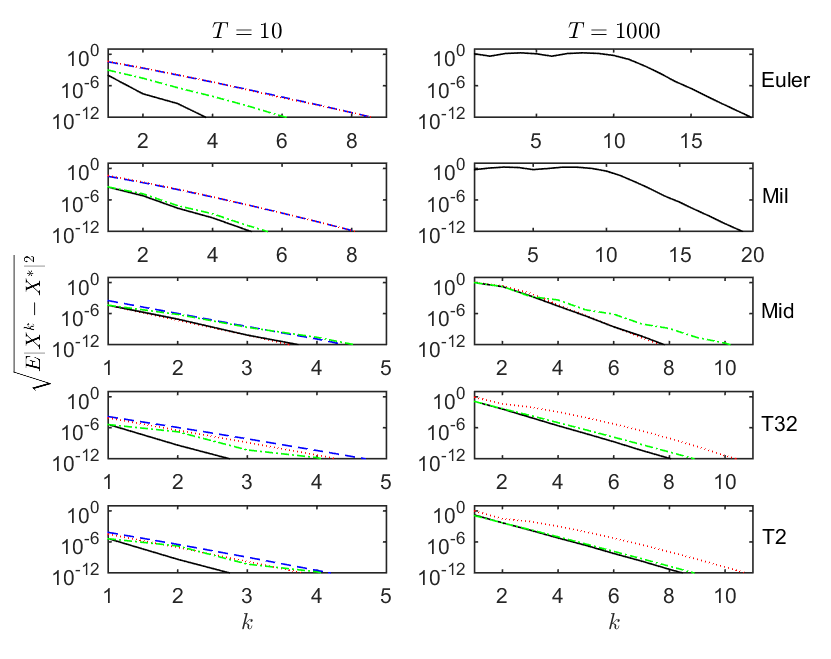}
	\caption{Kubo oscillator \eqref{kubo} with $c = 0.5$, $X_0=(1,0)$. Mean-square errors vs. iteration number $k$ for original parareal and parareal with projection algorithms using five propagators as $\F$ and $\G$ ($\Delta T = 0.1$, $J=100$). Left: $T=10$. Right: $T=1000$.}
	\label{fig;kubo}
\end{figure}

The convergence results of the parareal algorithms are shown in Figure \ref{fig;kubo}. The left part of this figure is the short time simulation (T=10), and the right part displays the long time simulation (T=1000). Each row of Figure \ref{fig;kubo} corresponds to a particular scheme which acts as the basic integrators ($\F$ and $\G$) in the parareal algorithm. In the case of short time test, we observe that all the schemes (with or without projection) converge properly, and the parareal algorithms with projection and using projection schemes as the $\G$ and $\F$ integrators, have the fastest convergence rate (the solid line). For the long time case, the common Euler and Milstein schemes without projection do not converge in the parareal algorithms, so we just plot the results of parareal with projection and EulerP or MilP in the first two rows of the right side of Figure \ref{fig;kubo}. However, the other high order methods still work in the corresponding parareal algorithms. Note that the Mid scheme preserves the quadratic conserved quantity \eqref{e;I_kubo}; thus, we do not need to use the MidP scheme in this test.

\begin{figure}[htbp]
	\includegraphics[width=.45\textwidth]{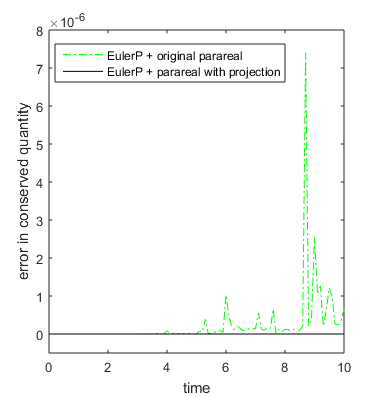}
	\includegraphics[width=.45\textwidth]{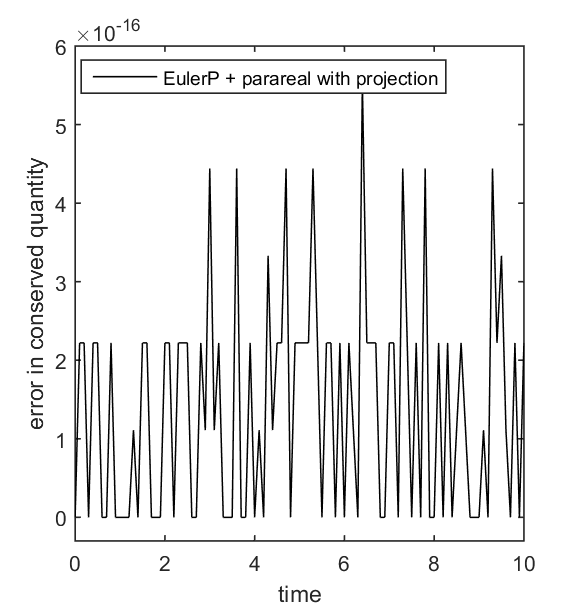}
	\caption{Errors in conserved quantity $I(x)$ \eqref{e;I_kubo} along numerical approximations by two kinds of parareal algorithms ($\F, \G$ = EulerP) after $k=2$ iterations. Left: the original parareal. Right: parareal with projection.}
	\label{fig;kubo_cq}
\end{figure}

From Figure \ref{fig;kubo}, it also turns out that with the help of projection method, the convergence rates for the original parareal and parareal with projection are similar if they both converge. To compare these two algorithms, we then estimate the errors of the conserved quantity $I(x)$ \eqref{e;I_kubo} in Figure \ref{fig;kubo_cq} where $T=10$. Here EulerP are chosen as the fine and coarse propagators $\F$ and $\G$. Other parameters are the same as those in Figure \ref{fig;kubo}. The left plot of Figure \ref{fig;kubo_cq} shows the errors of the original parareal and the parareal with projection after $k=2$ iterations, while the right solely demonstrates the later one. Therefore, although they both have good convergence property for SDEs with conserved quantity, the parareal with projection provides a much better reproduction of the preservation of the conserved quantity $I(x)$. In the case of other high mean-square order schemes, the results are just similar, so we omit them here.

\subsection{Stochastic pendulum}
Next, we restrict to a two-dimensional mathematical pendulum perturbed by two  multiplicative noises \cite{Cohen2014}
\begin{equation}\label{pend}
d
\begin{pmatrix}
X_1(t) \\
X_2(t) \\
\end{pmatrix}
=
\begin{pmatrix}
-\sin\big(X_2(t)\big) \\
X_1(t) \\
\end{pmatrix}
\Big( dt + c_1\circ dW_1(t) + c_2\circ dW_2(t)\Big),
\end{equation}
where $c_1$ and $c_2$ are real-valued parameters. It has a conserved quantity as follows
\begin{equation}\label{e;I_pend}
I(x,y) = \frac{1}{2}x^2 - \cos(y),
\end{equation}
which is a non-quadratic one unlike that of the first example.

\begin{figure}[htbp]
	\centering
	\includegraphics[width=.9\linewidth]{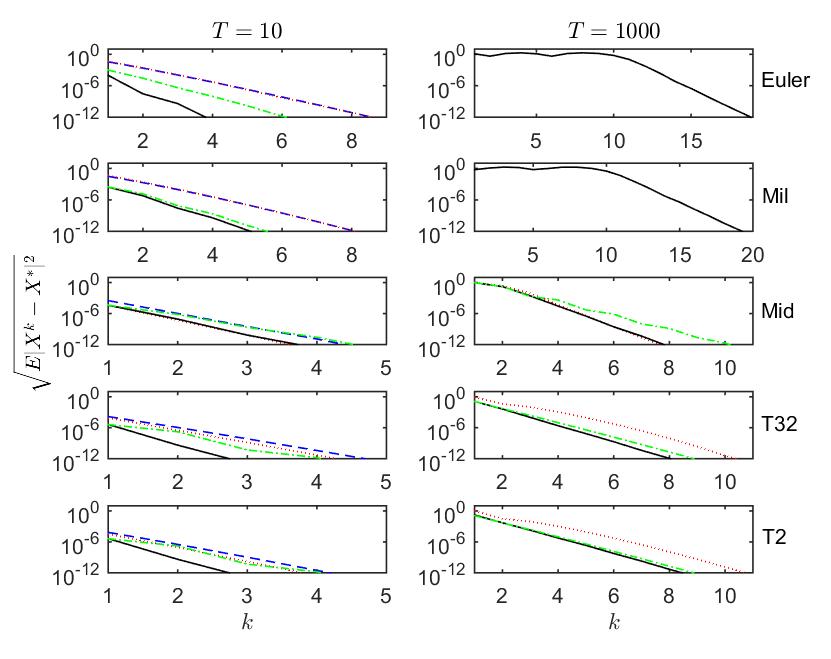}
	\caption{Stochastic pendulum with two multiplicative noises \eqref{pend} with $c_1 = 0.5$,  $c_2=0.1$, $X_0=(0.2,1)$. Mean-square errors vs. iteration number $k$ for original parareal and parareal with projection algorithms using five propagators as $\F$ and $\G$ ($\Delta T = 0.1$, $J=100$). Left: $T=10$. Right: $T=1000$.}
	\label{fig;pend}
\end{figure}

Setting $c_1 = 0.5, c_2 = 0.1$, and $X(0) = (0.2, 1)$, we thus get the results of convergence property in Figure \ref{fig;pend}. As before, the left and right part show the short time case and long time case, respectively. It is obvious that the results are similar to that of Figure \ref{fig;kubo}, except that Mid scheme can not preserve the conserved quantity \eqref{e;I_pend} (non-quadratic), so we consider the MidP scheme in the third row. In the case of $T=1000$, the original parareal algorithms without projection integrators are unable to reach proper error during the iteration process. Instead, the projection parareal algorithms using projection schemes as the $\G$ and $\F$ integrators converge properly. In addition, for Mid, T32P and T2 schemes with projection technique, the corresponding iteration numbers  are all less than 10.

\begin{figure}[htbp]
	\includegraphics[width=.45\textwidth]{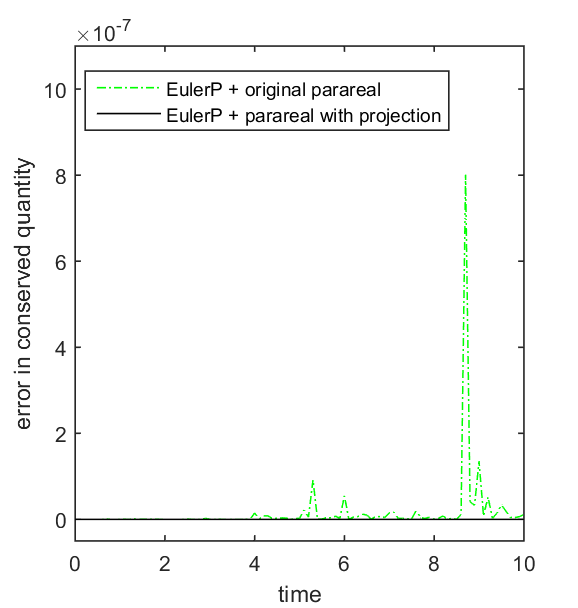}
	\includegraphics[width=.45\textwidth]{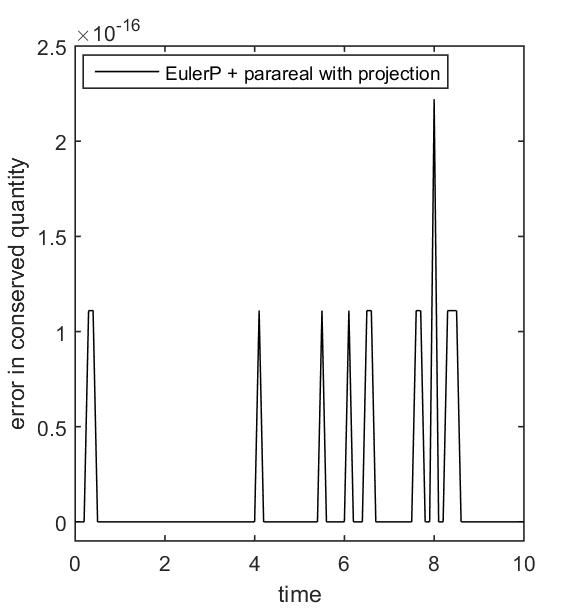}
	\caption{Errors in conserved quantity $I(x)$ \eqref{e;I_pend} along numerical approximations by two kinds of parareal algorithms ($\F, \G$ = EulerP) after $k=3$ iterations. Left: the original parareal. Right: parareal with projection}
	\label{fig;pend_cq}
\end{figure}

In addition, Figure \ref{fig;pend_cq} displays the errors in conserved quantity \eqref{e;I_pend} along the original parareal and parareal with projection,  where $T=10$, and other parameters are the same as those in the test of Figure \ref{fig;pend}. Here the fine and coarse propagators $\F, \G$ are all  EulerP. The left plot of Figure \ref{fig;pend_cq} shows the errors of the original parareal and the parareal with projection after $k=3$ iterations, while the right solely demonstrates the later one. Therefore, although they both have good convergence property for SDEs with conserved quantity, the parareal with projection provides a much better reproduction of the preservation of the conserved quantity \eqref{e;I_pend}.

\subsection{Stochastic cyclic Lotka-Volterra system}
Last we consider a three-dimensional cyclic Lotka-Volterra model
\begin{equation}\label{lk}
d\left(
\begin{array}{c}
X_1(t) \\
X_2(t) \\
X_3(t) \\
\end{array}
\right) =
\left(
\begin{array}{c}
X_1(t)\big(X_3(t)-X_2(t)\big) \\
X_2(t)\big(X_1(t)-X_3(t)\big) \\
X_3(t)\big(X_2(t)-X_1(t)\big) \\
\end{array}	\right) \Big(dt + c \circ dW(t)\Big),
\end{equation}
where $c$ is also a real-valued constant parameter. This system represents a chaotic environment consisting of three completing species \cite{Chen2014conservative,Hong2011}. And it is easy to check that system \eqref{lk} possesses two independent conserved quantities:
\begin{equation}\label{e;I_lk}
\begin{aligned}
I_1(x,y,z) &= x+y+z, \\
I_2(x,y,z) &= xyz.
\end{aligned}
\end{equation}

\begin{figure}[htbp]
	\centering
	\includegraphics[width=.9\linewidth]{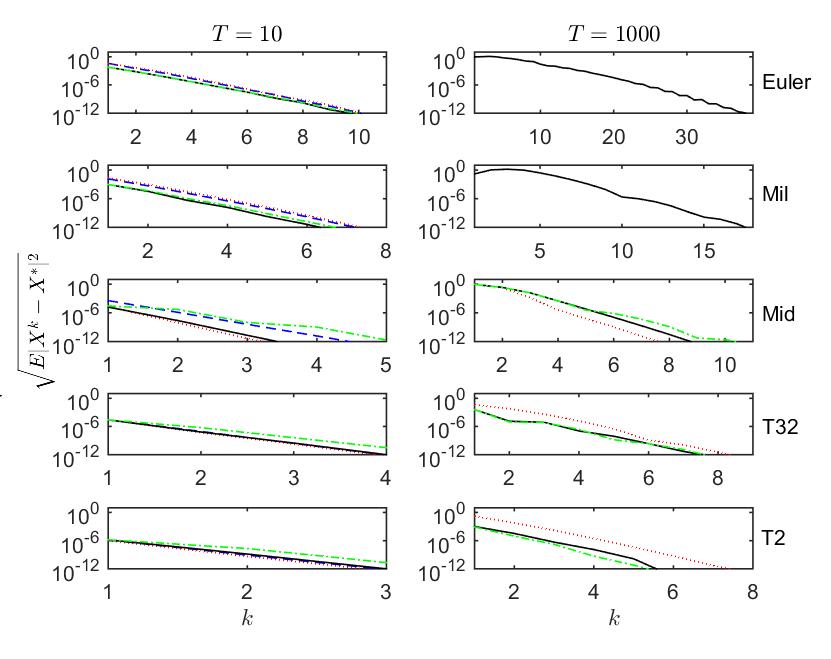}
	\caption{Stochastic cyclic Lotka-Volterra system \eqref{lk} with $c = 0.5$, $X_0=(1,2,1)$. Mean-square errors vs. iteration number $k$ for original parareal and parareal with projection algorithms using five propagators as $\F$ and $\G$ ($\Delta T = 0.01$, $J=100$). Left: $T=10$. Right: $T=1000$.}
	\label{fig;lk}
\end{figure}

By the conserved quantities above, the phase trajectory of the exact solution to \eqref{lk} is a closed curve in $\R^3$.
In this test, we choose parameter $c=0.5$ and initial value $X(0) = (1,2,1)$. In contrast to the previous two examples, we set $\Delta T = 0.01$ in order to investigate the long-term behavior of these methods. The convergence property for the corresponding parareal algorithms with different schemes are shown in Figure \ref{fig;lk}. Form the left part of it, we notice that all these algorithms are able to reach good convergence for $T=10$. However, figures in the right hand side show something different. First, for Euler and Mil schemes, we only plot the solid lines, i.e., projection performed in both propagators ($\F$ and $\G$) and correction step. Comparing these two figures, we observe that Mil type scheme can reach the stopping criterion \eqref{e;stopping} faster than the Euler one. That is to say, the Mil type scheme needs nearly one half of the iteration numbers compared to the Euler one, and achieve more accuracy. The last three figures show the convergence property of the Mid, T32 and T2 type schemes, respectively. Also, the original parareal algorithm with non-projection schemes are unable to meet the stopping criterion \eqref{e;stopping} during iteration, but if we use the projection schemes or parareal algorithm with projection, less than 10 iterations are needed for this example.

\begin{figure}[htbp]
	\includegraphics[width=.45\textwidth]{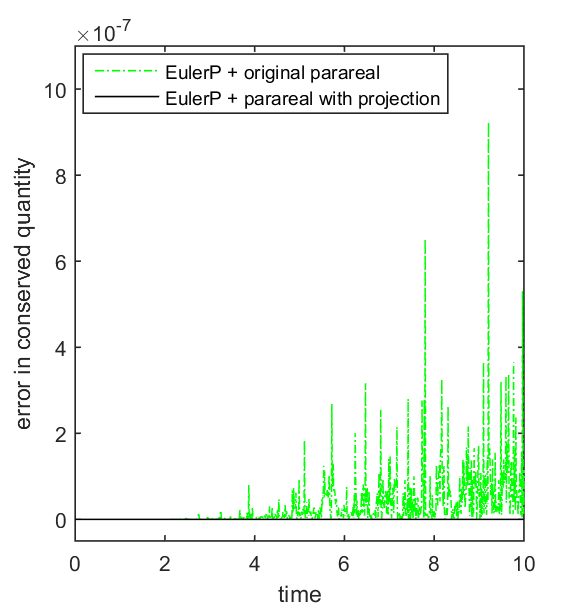}
	\includegraphics[width=.45\textwidth]{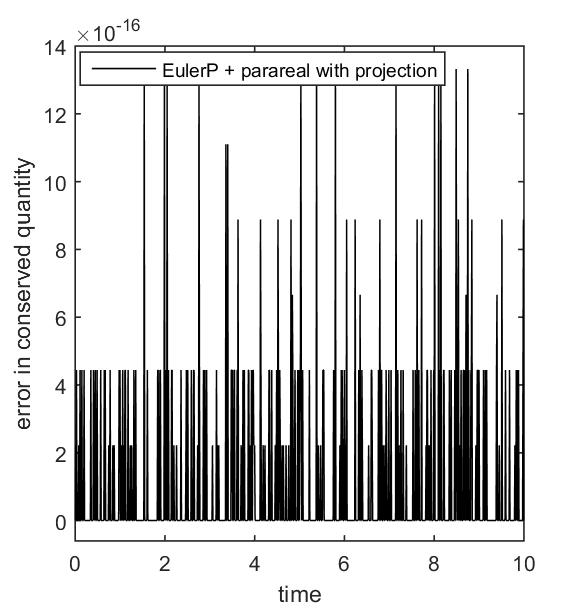}
	\caption{Errors in conserved quantity $I(x)$ \eqref{e;I_lk} along numerical approximations by two kinds of parareal algorithms ($\F, \G$ = EulerP) after $k=5$ iterations. Left: the original parareal. Right: parareal with projection}
	\label{fig;lk_cq}
\end{figure}

Errors in the two conserved quantities \eqref{e;I_lk} along the original parareal and parareal with projection algorithms are shown in Figure \ref{fig;lk_cq}, where $T=10$, and other parameters are the same as those in the test of Figure \ref{fig;lk}. Also the $\F$ and $\G$ all choose the EulerP scheme same as the previous examples. The left plot of Figure \ref{fig;lk_cq} shows the errors of the original parareal and the parareal with projection after $k=5$ iterations, while the right solely demonstrates the later one. Therefore, although they both have good convergence property for SDEs with conserved quantity, the parareal with projection provides a much better reproduction of the preservation of the conserved quantity \eqref{e;I_lk}.

\section{Conclusion}
In conclusion, we investigate the possibility of applying parallel-in-time technique to SDEs with conserved quantities by combing the projection methods and a version of the parareal algorithm. For this kind of system, projection methods can be used to guarantee that the numerical approximations preserve certain conserved quantities exactly. However, the long-time simulation of this problem is still challenging, and it is inevitable to take the parallel algorithms into consideration. With the help of the parareal algorithm with projection, we obtain an effective parallel-in-time approach which maintains the geometric property to simulate SDEs with conserved quantities. In the numerical experiments, three systems, linear or non-linear, are performed by parareal algorithm with and without projection technique, respectively. From the numerical results, we can conclude that the parareal algorithm is of fast convergence with few iterations, and with the help of projection method it shows advantages in preserving conserved quantities. This paper mainly focuses on the numerical simulation of this efficient parareal algorithm for SDEs with conserved quantities. However, there are lack of corresponding theoretical analysis of this algorithm. Thus, we will continue to study the convergence property and other numerical behaviors in the further works.
%

\bibliographystyle{elsarticle-num}
\bibliography{Parareal}

\begin{thebibliography}{10}
\expandafter\ifx\csname url\endcsname\relax
  \def\url#1{\texttt{#1}}\fi
\expandafter\ifx\csname urlprefix\endcsname\relax\def\urlprefix{URL }\fi
\expandafter\ifx\csname href\endcsname\relax
  \def\href#1#2{#2} \def\path#1{#1}\fi

\bibitem{Lions2001}
J.-L. Lions, Y.~Maday, G.~Turinici, {A ``parareal" in time discretization of
  PDE's}, Comptes Rendus l'Acad{\'{e}}mie des Sci. - Ser. I - Math. 332~(7)
  (2001) 661--668.
\newblock \href {http://dx.doi.org/10.1016/S0764-4442(00)01793-6}
  {\path{doi:10.1016/S0764-4442(00)01793-6}}.

\bibitem{Bal2002}
G.~Bal, Y.~Maday, {A "Parareal" Time Discretization for Non-Linear PDE's with
  Application to the Pricing of an American Put}, in: Recent Dev. domain
  Decompos. methods, Springer Berlin Heidelberg, 2002, pp. 189--202.
\newblock \href {http://dx.doi.org/10.1007/978-3-642-56118-4_12}
  {\path{doi:10.1007/978-3-642-56118-4_12}}.

\bibitem{Wu2015}
S.~Wu, T.~Zhou, {Convergence Analysis for Three Parareal Solvers}, SIAM J. Sci.
  Comput. 37~(2) (2015) A970--A992.
\newblock \href {http://dx.doi.org/10.1137/140970756}
  {\path{doi:10.1137/140970756}}.

\bibitem{Maday2002}
Y.~Maday, G.~Turinici, A parareal in time procedure for the control of partial
  differential equations, Comptes Rendus Mathematique 335~(4) (2002) 387--392.
\newblock \href {http://dx.doi.org/10.1016/s1631-073x(02)02467-6}
  {\path{doi:10.1016/s1631-073x(02)02467-6}}.

\bibitem{2005Barth}
T.~J. Barth, M.~Griebel, D.~E. Keyes, R.~M. Nieminen, D.~Roose, T.~Schlick,
  R.~Kornhuber, R.~Hoppe, J.~P{\'{e}}riaux, O.~Pironneau, O.~Widlund, J.~Xu
  (Eds.), Domain Decomposition Methods in Science and Engineering, Springer
  Berlin Heidelberg, 2005.
\newblock \href {http://dx.doi.org/10.1007/b138136}
  {\path{doi:10.1007/b138136}}.

\bibitem{Dai2013}
X.~Dai, C.~{Le Bris}, F.~Legoll, Y.~Maday, {Symmetric parareal algorithms for
  Hamiltonian systems}, ESAIM Math. Model. Numer. Anal. 47~(3) (2013) 717--742.
\newblock \href {http://dx.doi.org/10.1051/m2an/2012046}
  {\path{doi:10.1051/m2an/2012046}}.

\bibitem{Gander2014}
M.~J. Gander, E.~Hairer, {Analysis for parareal algorithms applied to
  Hamiltonian differential equations}, J. Comput. Appl. Math. 259 (2014) 2--13.
\newblock \href {http://dx.doi.org/10.1016/j.cam.2013.01.011}
  {\path{doi:10.1016/j.cam.2013.01.011}}.

\bibitem{Gardiner2009}
C.~Gardiner, {Stochastic Methods: A Handbook for the Natural and Social
  Sciences}, 4th Edition, Springer-Verlag Berlin Heidelberg, 2009.

\bibitem{Oksendal2003}
B.~{\O}ksendal, {Stochastic Differential Equations}, Springer-Verlag, Berlin,
  2003.
\newblock \href {http://dx.doi.org/10.1007/978-3-642-14394-6}
  {\path{doi:10.1007/978-3-642-14394-6}}.

\bibitem{Burrage2004}
K.~Burrage, P.~M. Burrage, T.~Tian, {Numerical methods for strong solutions of
  stochastic differential equations: an overview}, Proc. R. Soc. A Math. Phys.
  Eng. Sci. 460~(2041) (2004) 373--402.
\newblock \href {http://dx.doi.org/10.1098/rspa.2003.1247}
  {\path{doi:10.1098/rspa.2003.1247}}.

\bibitem{Higham2001}
D.~J. Higham, {An algorithmic introduction to numerical simulation of
  stochastic differential equations}, SIAM Rev. 43~(3) (2001) 525--546.
\newblock \href {http://dx.doi.org/http://dx.doi.org/10.1137/S0036144500378302}
  {\path{doi:http://dx.doi.org/10.1137/S0036144500378302}}.

\bibitem{Kloeden1992}
P.~E. Kloeden, E.~Platen, {Numerical Solution of Stochastic Differential
  Equations}, Springer-Verlag, Berlin, 1992.
\newblock \href {http://dx.doi.org/10.1007/978-3-662-12616-5}
  {\path{doi:10.1007/978-3-662-12616-5}}.

\bibitem{Hairer2006}
E.~Hairer, C.~Lubich, G.~Wanner, {Geometric Numerical Integration:
  Structure-Preserving Algorithms for Ordinary Differential Equations},
  Springer, Berlin, 2002.
\newblock \href {http://dx.doi.org/10.1007/3-540-30666-8}
  {\path{doi:10.1007/3-540-30666-8}}.

\bibitem{Milstein2004}
G.~Milstein, M.~Tretyakov, {Stochastic Numerics for Mathematical Physics},
  Springer-Verlag, Berlin, 2004.
\newblock \href {http://dx.doi.org/10.1007/978-3-662-10063-9}
  {\path{doi:10.1007/978-3-662-10063-9}}.

\bibitem{Bal2003}
G.~Bal, {Parallelization in time of (stochastic) ordinary differential
  equations}, Preprint (2003) 1--23.

\bibitem{Engblom2009}
S.~Engblom, {Parallel in Time Simulation of Multiscale Stochastic Chemical
  Kinetics}, Multiscale Model. Simul. 8~(1) (2009) 46--68.
\newblock \href {http://dx.doi.org/10.1137/080733723}
  {\path{doi:10.1137/080733723}}.

\bibitem{Dai2013a}
X.~Dai, Y.~Maday, {Stable Parareal in Time Method for First-and Second-Order
  Hyperbolic Systems}, SIAM J. Sci. Comput. 35~(1) (2013) A52--A78.
\newblock \href {http://dx.doi.org/10.1137/110861002}
  {\path{doi:10.1137/110861002}}.

\bibitem{Wu2011}
S.~Wu, Z.~Wang, C.~Huang, {Analysis of mean-square stability of the parareal
  algorithm}, Math. Numer. Sinica 33~(2) (2011) 113--124.

\bibitem{Zhou2016}
W.~Zhou, L.~Zhang, J.~Hong, S.~Song,
  \href{http://link.springer.com/10.1007/s10543-016-0614-0}{{Projection methods
  for stochastic differential equations with conserved quantities}}, BIT 56~(4)
  (2016) 1497--1518.
\newblock \href {http://arxiv.org/abs/arXiv:1601.04157v1}
  {\path{arXiv:arXiv:1601.04157v1}}, \href
  {http://dx.doi.org/10.1007/s10543-016-0614-0}
  {\path{doi:10.1007/s10543-016-0614-0}}.
\newline\urlprefix\url{http://link.springer.com/10.1007/s10543-016-0614-0}

\bibitem{Cohen2014}
D.~Cohen, G.~Dujardin, {Energy-preserving integrators for stochastic Poisson
  systems}, Commun. Math. Sci. 12~(8) (2014) 1523--1539.
\newblock \href {http://dx.doi.org/10.4310/CMS.2014.v12.n8.a7}
  {\path{doi:10.4310/CMS.2014.v12.n8.a7}}.

\bibitem{Milstein2002mul}
G.~N. Milstein, Y.~M. Repin, M.~V. Tretyakov, {Numerical methods for stochastic
  systems preserving symplectic structure}, SIAM J. Numer. Anal. 40~(4) (2002)
  1583--1604.
\newblock \href {http://dx.doi.org/10.1137/S0036142901395588}
  {\path{doi:10.1137/S0036142901395588}}.

\bibitem{Chen2014conservative}
C.~Chen, D.~Cohen, J.~Hong, {Conservative methods for stochastic differential
  equations with a conserved quantity}, Int. J. Numer. Anal. Model. 13~(3)
  (2016) 435--456.

\bibitem{Hong2011}
J.~Hong, S.~Zhai, J.~Zhang, {Discrete gradient approach to stochastic
  differential equations with a conserved quantity}, SIAM J. Numer. Anal.
  49~(5) (2011) 2017--2038.
\newblock \href {http://dx.doi.org/10.1137/090771880}
  {\path{doi:10.1137/090771880}}.

\end{thebibliography}

\end{document}